\documentclass[letterpaper,10pt]{article}

\usepackage[latin1]{inputenc}
\usepackage{amsfonts}
\usepackage{amsmath}
\usepackage{amssymb}
\usepackage{amsthm}
\usepackage[american]{babel}
\usepackage[dvips]{psfrag,graphics}
\usepackage{enumerate}

\newtheorem{theorem}{Theorem}[section]
\newtheorem*{numarasizteo}{Theorem}
\newtheorem{proposition}[theorem]{Proposition}
\newtheorem{definition}{Definition}
\newtheorem{lemma}[theorem]{Lemma}
\newtheorem{corollary}[theorem]{Corollary}

\newenvironment{problem}[1][Problem]{\textbf{#1.} }

\newenvironment{listi}{\begin{enumerate}[\upshape(i)]\setlength{\itemsep}{1pt}\setlength{\parskip}{0pt}\setlength{\parsep}{0pt}}
{\end{enumerate}}

\newenvironment{lista}{\begin{enumerate}[\upshape(a)]\setlength{\itemsep}{1pt}\setlength{\parskip}{0pt}\setlength{\parsep}{0pt}}
{\end{enumerate}}

\newcommand{\Ha}{\mathcal{H}}
\newcommand{\la}{\lambda}
\newcommand{\ga}{\gamma}
\newcommand{\al}{\alpha}

\newcommand{\ep}{\epsilon}
\newcommand{\de}{\delta}
\newcommand{\om}{\omega}
\newcommand{\Om}{\Omega}

\newcommand{\ip}{\mathcal{IP}}

\title{On the arithmetic sums of Cantor sets}
\author{Kemal Ilgar Ero\u{g}lu}

\begin{document}
\maketitle
\begin{abstract}
Let $C_\la$ and $C_\ga$ be two affine Cantor sets in $\mathbb{R}$ with similarity dimensions $d_\la$ and $d_\ga$, respectively. We define an analog of the Bandt-Graf condition for self-similar systems and use it to give necessary and sufficient conditions for having $\Ha^{d_\la+d_\ga}(C_\la + C_\ga)>0$ where $C_\la + C_\ga$ denotes the arithmetic sum of the sets. We use this result to analyze the orthogonal projection properties of sets of the form $C_\la \times C_\ga$. We prove that for Lebesgue almost all directions $\theta$ for which the projection is not one-to-one, the projection has zero $(d_\la + d_\ga)$-dimensional Hausdorff measure. We demonstrate the results on the case when $C_\la$ and $C_\ga$ are the middle-$(1-2\la)$ and middle-$(1-2\ga)$ sets.
\end{abstract}

\section{Introduction}

The work presented in this paper is inspired by several results which have connection to the study of certain dynamical systems. In the work of J.~Palis \cite{palis}, arithmetic sums of Cantor sets came up in the study of homoclinic tangencies. The main question was whether such sums contained intervals, provided they have positive Lebesgue measure. This lead to similar questions involving affine Cantor sets in the line. Y.~Peres and B.~Solomyak \cite{SolInt} obtained a result on the dimensions of arithmetic sums under certain conditions. It is then natural to investigate the Hausdorff measure of these sums in their dimensions. This paper aims to answer some questions in this direction.

We develop some criteria that serve as adaptations of the Bandt-Graf condition for self-similar systems to the case of arithmetic sums of affine Cantor sets. We define the concepts of ``$\ep$-square'' and ``$\ep$-relative closeness''. We use a ``doubling'' argument (see Lemma~\ref{cl4}) inspired by the techniques in \cite{SolNonlinear} that has also been used in \cite{irrat_rot} in a similar manner. We formulate our main result as follows:

\begin{problem}[Theorem \ref{Main1}] \emph{Let $d_\la,d_\ga$ be the similarity dimensions of the affine Cantor sets $C_\la,C_\ga \subset\mathbb{R}$, respectively. Then, $\mathcal{H}^{d_\gamma + d_\lambda}(C_\la + C_\ga)=0$ if and only if for any $\epsilon > 0$ there exist two distinct $\epsilon$-relatively close $\ep$-squares.}
\end{problem}

After proving this result, we give some applications on the orthogonal projections of sets of type $C_\la \times C_\ga$. Orthogonal projections $\Pi_\theta(C_\la \times C_\ga)$ onto lines can be identified (up to rotations and scaling) with arithmetic sums of the form $C_\la + tC_\la$. In particular, we prove that for Lebesgue almost all directions $\theta$ for which the projection is not one-to-one, we have $\mathcal{H}^{d_\gamma + d_\lambda}(\Pi_\theta(C_\la + C_\ga))=0$. Also, we prove that if $C_\la$ and $C_\ga$ satisfy the strong separation condition and $d_\la+d_\ga < 1/2$, then for Lebesgue almost all projections we have $\mathcal{H}^{d_\gamma + d_\lambda}(\Pi_\theta(C_\la + C_\ga))>0$.

Finally, we consider the concrete example where $C_\la$ and $C_\ga$ are the middle-$(1-2\la)$ and middle-$(1-2\ga)$ sets. The case when $d_\la + d_\ga > 1 $ has already been studied in detail by several authors. Our work provides new information on the case $d_\la + d_\ga < 1 $.
\section{Arithmetic sums of affine Cantor sets}\label{prob1}
\subsection{Background}

By an \emph{affine Cantor} set (in $\mathbb{R}$) we mean the attractor of an iterated function system whose maps are linear contractions. We say the set is homogeneous if all the contraction rates are the same. A \emph{regular Cantor set} is the attractor of an iterated function system with smooth contractions.

In \cite{PalisTakens}, J.~Palis and F.~Takens asked whether it was true that arithmetic sums of regular Cantor sets generically, if not always, contained intervals if they had positive Lebesgue measure. They also asked the same question for affine Cantor sets.

The question about regular sets was answered by A.~Moreira and J.~Yoccoz in \cite{Moreira} (they proved that the statement is generically true). Sannami \cite{Sannami} constructed an affine Cantor set whose self-difference is a Cantor set with positive Lebesgue measure. The question about the generic sums of affine Cantor sets is still an open problem.

A result that was related to this question was obtained by Y.~Peres and B.~Solomyak in \cite{SolInt}. They considered homogeneous Cantor set families of the form 
\[
\left\{C_\la^\mathcal{D(\la)}\right\} = \left\{ \sum_{n=0}^\infty d_{i_n}(\la) \la^n \mid i_n\in\{1,\ldots,m\} \right\}
\]
where \mbox{$\mathcal{D}(\la)=\{d_1 (\la), \ldots,d_m(\la)\}$} is such that $d_i (\la) \in C^1 [0,1]$ for $i=1,\ldots,m$. They obtained the following result:

\begin{numarasizteo}[Peres, Solomyak]\label{SolInt1}
Suppose $K$ is a compact set on the real line and $J \subset (0,1)$ is an interval such that the family $\{C_\la^\mathcal{D(\la)}\}$ satisfies the strong separation condition for all $\la \in J$. Then
 
 \begin{lista}
 \item for almost every $\la \in J$ such that $\dim_H C_\la^\mathcal{D(\la)} + \dim_H K < 1$, we have 
 \[
 \dim_H \left(C_\la^\mathcal{D(\la)} + K\right) = \dim_H C_\la^\mathcal{D(\la)} + \dim_H K.
 \]
 \item for almost every $\la \in J$ such that $\dim_H C_\la^\mathcal{D(\la)} + \dim_H K >1$, the set $C_\la^\mathcal{D(\la)} + K$ has positive Lebesgue measure.
 \end{lista}
\end{numarasizteo}

Part (a) implies that if the Hausdorff dimensions of two homogeneous Cantor sets with strong separation add up to a number less than one, then, typically, the dimension of their sum is the sum of their dimensions. A question one can then ask is about the measure of these arithmetic sums, in the appropriate dimension. More precisely, one can ask when these sums have positive measure in the ``right'' dimension.

For self-similar sets, there are several equivalent ways to characterize those sets that have positive Hausdorff measure in their similarity dimension, such as the Open Set Condition (see \cite{Fal2}) or the Bandt-Graf condition \cite{BG}. In this section we try to carry the ideas of the latter method over to the case of arithmetic sums and find a necessary and sufficient condition for having positive Hausdorff measure in an appropriate dimension (namely, the sum of similarity dimensions).

\subsection{The result}

We consider two self-similar sets $C_\la$ and $C_\ga$ in the real line that are the attractors of iterated function systems whose maps are $\{F_1,\ldots,F_A\}$ and $\{G_1,\ldots,G_B\}$, respectively. For each $i=1,\ldots,A$ and $j=1,\ldots,B$ write 
\begin{gather*} 
F_i(x)=O^\la_i \la_i x + b_i, \\
G_j(x)=O^\ga_j \ga_j x + c_j.
\end{gather*}
Here the $O^\la_i$ and $O^\ga_j$ are $+1$ or $-1$ and $1>\la_i,\ga_j >0$. Define
\begin{equation}\label{def_rmin}
r_{min}=\min\{\la_1,\ldots,\la_A,\ga_1,\ldots,\ga_B\}.
\end{equation} 

Let $\Omega_\la=\{1,2,\ldots,A\}^\mathbb{N}$, $\Omega_\ga=\{1,2,\ldots,B\}^\mathbb{N}$ be spaces of infinite sequences of $A$ and $B$ digits. Then there are natural projections $\Pi_\la,\Pi_\ga$ from these spaces onto the sets $C_\la$ and $C_\ga$. For example for $\Pi_\la$ we can take
\[
\Pi_\la:\ \om=(\om_1 \om_2 \om_3\ldots)\ \longrightarrow\ \lim_{n\to\infty} F_{\om_1} \circ F_{\om_2} \circ F_{\om_3} \cdots F_{\om_n}(0)
\]
(see \cite{Fal2} for details). If the sequence $\om$ is mapped to $x$, we say $\om$ is an \emph{address} of $x$. Finite sequences of digits will be called \emph{words}. A bar above a digit or a word will mean that the word or digit repeats infinitely; for example $0\bar{1}$ will mean $0111\ldots$ 

If $u$ is a word of the alphabet used for, say, $\Om_\la$, then by $[u]$ we will denote the set of all elements in $\Om_\la$ starting with the word $u$ (such sets are called \emph{cylinder sets}). Given a word $u=u_1 \cdots u_n$, write
\[
F_u = F_{u_1} \circ F_{u_2} \circ \cdots \circ F_{u_n}.
\]
We give a similar definition for $G_u$. We also define $O^\la_u=O^\la_{u_1} \cdots O^\la_{u_n}$ (again, similar definition for $O^\ga_u$). We will denote $\la_u=\la_{u_1}\la_{u_2}\cdots\la_{u_n}$ and $\ga_u=\ga_{u_1}\ga_{u_2}\cdots\ga_{u_n}$. We will use $|u|$ to denote the length of the word $u$.

Let $D_\la$ and $D_\ga$ be the diameters of $C_\la$ and $C_\ga$. Write
\[
D= \max \{D_\la,D_\ga\}.
\]

Let $ \Omega= \Omega_\la \times \Omega_\ga$ and let $\Pi=\Pi_\la \times \Pi_\ga$ be the ``projection'' from $\Om$ onto $C_\la \times C_\ga \subset \mathbb{R}^2$. Let $L$ be the map from $\Om$ on $C_\la + C_\ga$ given by 

\begin{equation}\label{def_L}
L(\om,\tau) = \Pi_\la(\om) + \Pi_\ga(\tau).
\end{equation}
We observe the simple geometric fact that the orthogonal projection of $C_\la \times C_\ga$ onto the line $y=x$ in the plane is the same as $C_\la + C_\ga$, up to scaling by a factor of $\sqrt{2}$. Therefore one can also visualize the map $L$ as the composition of $\Pi$ with the orthogonal projection onto the line $y=x$.

Finally, we will write $[u \times v]$ to denote $[u] \times [v]$ where $u$ and $v$ are words from the alphabets of $\Om_\la$ and $\Om_\ga$, respectively. We will refer to $u \times v$ as a word (of $\Om$), corresponding to the ``cylinder set'' $[u \times v]$.

\bigskip

We give two definitions before stating the main result:

\begin{definition}
A word $u \times v$ of $\Om$ is called a $\de$-square if $\la_u / \ga_v \in (e^{-\de},e^\de)$.
\end{definition}

\begin{definition}\label{def2}
We say $u_1 \times v_1$ and $u_2 \times v_2$ are $\epsilon$-relatively close if 

\begin{listi}

\item $\la_{u_1}/\la_{u_2} \in (e^{-\ep},e^\ep)$ and $\ga_{v_1}/\ga_{v_2} \in (e^{-\ep},e^\ep)$;

\item $O^\la_{u_1} =O^\la_{u_2}$ and $O^\ga_{v_1} =O^\ga_{v_2}$;
\item $|L(u_{1} \bar{1},v_{1} \bar{1}) - L(u_{2} \bar{1},v_{2} \bar{1})| < \epsilon D \min \{\la_{u_1}, \la_{u_2}, \ga_{v_1}, \ga_{v_2}\}.$
\end{listi}
\end{definition} 

\bigskip
Recall that the similarity dimension of $C_\la$ is the real number $d_\la$ such that $\sum_{i=1}^A \la_i^{d_\la}=1$. A similar definition is given for $d_\ga$.

The main result is as follows:
\begin{theorem}\label{Main1}
Let $d_\la,d_\ga$ be the similarity dimensions of $C_\la,C_\ga$, respectively. Then, $\mathcal{H}^{d_\gamma + d_\lambda}(C_\la + C_\ga)=0$ if and only if for any $\epsilon > 0$ there exist two distinct $\epsilon$-relatively close $\ep$-squares.
\end{theorem}

The general strategy of the proof is similar to the approach in \cite{SolNonlinear}. We start with a sequence of preliminary results.

\begin{lemma}\label{cl1}
Let $u_1 \times v_1$ and $u_2 \times v_2$ be $\epsilon$-relatively close $\delta$-squares. Let $s \times t$ be a $\delta^\prime$-square with $O^\la_s=O^\ga_t$. Then the words  $su_1 \times tv_1$ and $su_2 \times tv_2$ are \mbox{$( \epsilon e^{\delta^\prime} + 2  |e^{\delta^\prime} -1 | e^{\de^\prime +\de+\ep} \la_{u_1}^{-1} )$} -relatively close $(\delta + \delta^\prime)$-squares.	
\end{lemma}
\begin{proof}
Let $(a,b)=\Pi(u_{1}\bar{1},v_{1} \bar{1})$, $(c,d)= \Pi(u_{2} \bar{1}, v_{2} \bar{1})$. By translating if necessary, we may assume $C_\ga$ is contained in $[0,D]$ so that $|b|,|d|\leq D$. Then
\[
\begin{split}
&|L(su_{1} \bar{1}, tv_{1} \bar{1}) - L(su_{2} \bar{1}, tv_2\bar{1})|\\
&=|\Pi_\la(su_{1} \bar{1}) +\Pi_\ga(tv_{1} \bar{1})- \Pi_\la(su_{2} \bar{1})  - \Pi_\ga( tv_2\bar{1})|\\
&= |O^\la_s \la_{s}a + O^\ga_t\ga_{t}b -(O^\la_s\la_{s}c + O^\ga_t\ga_{t}d)| \\
&= |\la_{s}a + \ga_{t}b -(\la_{s}c + \ga_{t}d)| \\
& \leq \la_{s}|a+b - (c+d)| + 
|b|\la_{s}\left|1-\frac{\ga_{t}}{\la_{s}}\right| + |d|\la_{s}\left|1-\frac{\ga_{t}}{\la_{s}}\right| \\
& \leq \la_{s} \ep D \min\{\la_{u_1},\la_{u_2},\ga_{v_1},\ga_{v_2}\} + 2 D \la_{s}|e^{\de^\prime}-1|
\end{split}
\]
Observing that
\begin{equation}\label{relations}
\frac{\la_s}{\ga_t}\in(e^{-\delta^\prime},e^{\delta^\prime}), \quad \frac{\la_{u_i}}{\ga_{v_i}}\in(e^{-\delta},e^{\delta}),\,i=1,2 \ \ \text{  and  }\ \ 
\frac{\la_{u_1}}{\la_{u_2}},\frac{\ga_{v_1}}{\ga_{v_2}}\in(e^{-\ep},e^{\ep}) \quad
 \end{equation}
we obtain

\[
\begin{split}
&|L(su_{1} \bar{1}, tv_{1} \bar{1}) - L(su_{2} \bar{1}, tv_2\bar{1})|\\
& \leq \epsilon e^{\delta^\prime} D \min\{\la_{su_1},\la_{su_2},\ga_{tv_1},\ga_{tv_2}\} + 2 D \la_{s}| e^{\de^\prime}-1| e^{\delta + \ep} \la_{u_1}^{-1} \min\{\la_{u_1},\la_{u_2},\ga_{v_1},\ga_{v_2}\} \\
&\leq \left( \epsilon e^{\delta^\prime} + 2  |e^{\delta^\prime} -1 | e^{\de^\prime +\de+\ep} \la_{u_1}^{-1} \right) D \min\{\la_{su_1},\la_{su_2},\ga_{tv_1},\ga_{tv_2}\}.
\end{split}
\]

Finally, one can easily see that $\la_{su_i}/\ga_{tv_i} \in (e^{-(\delta+\delta^\prime)},e^{\delta + \delta^\prime})$, $i=1,2$ and also that $O^\la_{su_1}=O^\la_{su_2}$, $O^\ga_{tv_1}=O^\ga_{tv_2}$.
\end{proof}

\begin{lemma}\label{cl2}
Let $u_1 \times v_1$ and $u_2 \times v_2$ be $\ep$-relatively close $\de$-squares. Let $s\times t$ be a $\delta^\prime$-square such that $O^\la_s=O^\ga_t$. Then $u_{1}s \times v_{1}t$ and $u_{2}s \times v_{2}t$ are \mbox{$e^{\de^\prime} \la_s^{-1} (\ep + 2 e^{\ep +\delta} |e^\ep-1|)$}-relatively close $(\delta + \delta^\prime)$-squares.
\end{lemma}

\begin{proof}
It is clear that $u_{1}s \times v_{1}t$ and $u_{2}s \times v_{2}t$ are $(\delta + \delta^\prime)$-squares. Let $c=\Pi_\la(s\bar{1})-\Pi_\la(\bar{1})$, $d=\Pi_\ga(t\bar{1})-\Pi_\ga(\bar{1})$. Observe that

\begin{gather*}
\Pi_\la(u_1 s\bar{1})-\Pi_\la(u_1\bar{1}) = O^\la_{u_1} \la_{u_1} c, \\
\Pi_\la(u_2 s\bar{1})-\Pi_\la(u_2\bar{1}) = O^\la_{u_2} \la_{u_2} c, \\
\Pi_\ga(v_1 t\bar{1})-\Pi_\ga(v_1\bar{1}) = O^\ga_{v_1} \ga_{v_1} d, \\
\Pi_\ga(v_2 t\bar{1})-\Pi_\ga(v_2\bar{1}) = O^\ga_{v_2} \ga_{v_2} d.
\end{gather*}
Using these equalities and the definition of $L$, we get
\[
\begin{split}
&|L(u_{1}s \bar{1}, v_{1}t \bar{1}) - L(u_{2}s \bar{1}, v_2t\bar{1})|\\
&=\left|L(u_{1}\bar{1}, v_{1}\bar{1}) - L(u_{2}\bar{1}, v_2\bar{1}) + O^\la_{u_1}\la_{u_1}c-O^\la_{u_2}\la_{u_2}c + O^\ga_{v_1}\ga_{v_1}d-O^\ga_{v_2}\ga_{v_2}d\right|\\
&\leq \ep D \min\{\la_{u_1},\la_{u_2},\ga_{v_1},\ga_{v_2}\} + |c||\la_{u_1} - \la_{u_2}| + |d||\ga_{v_1}-\ga_{v_2}|\\
& \leq \ep D \min\{\la_{u_1},\la_{u_2},\ga_{v_1},\ga_{v_2}\} + D\la_{u_1}\left|1-\frac{\la_{u_1}}{\la_{u_2}}\right| + D\ga_{v_1}\left|1 -\frac{\ga_{v_1}}{\ga_{v_2}}\right| \\
& \leq \ep D \min\{\la_{u_1},\la_{u_2},\ga_{v_1},\ga_{v_2}\} + D\la_{u_1 }|e^\ep-1| + D\ga_{v_1}|e^\ep-1|
\end{split}
\]
Again by using (\ref{relations}) we can write
\[
\begin{split}
|L(u_{1}s \bar{1}, v_{1}t \bar{1}) &- L(u_{2}s \bar{1}, v_2t\bar{1})|\\
&\leq \ep e^{\delta^\prime} \la_s^{-1} D \min\{\la_{u_1 s},\la_{u_2 s},\ga_{v_1 t},\ga_{v_2 t}\}\\
&\quad + 2 D e^{\de + \de^\prime + \ep} |e^\ep-1| \la_s^{-1}\min\{\la_{u_1 s},\la_{u_2 s},\ga_{v_1 t},\ga_{v_2 t}\}  \\
& \leq e^{\de^\prime} \la_s^{-1} (\ep + 2 e^{\ep +\delta} |e^\ep-1|) D \min\{\la_{u_1 s},\la_{u_2 s},\ga_{v_1 t},\ga_{v_2 t}\}.
\end{split}
\]

Finally, checking that $O^\la_{u_1 s} = O^\la_{u_2 s}$ and $O^\ga_{v_1 t} = O^\ga_{v_2 t}$ are trivial.
\end{proof}

\begin{lemma}\label{cl3}
Let $0<\ep<\log2$. If $u_1 \times v_1$ and $u_2 \times v_2$ are $\epsilon$-relatively close, and if $u_2 \times v_2$ and $u_3 \times v_3$ are $\epsilon$-relatively close, then $u_1 \times v_1$ and $u_3 \times v_3$ are $4\epsilon$-relatively close.
\end{lemma}
\begin{proof}
\[
\begin{split}
&|L((u_{3} \bar{1}, v_{3} \bar{1})) - L((u_{1} \bar{1}, v_1\bar{1}))|\\ &\leq |L((u_{3} \bar{1}, v_{3} \bar{1})) - L((u_{2} \bar{1}, v_2\bar{1}))| + |L((u_{2} \bar{1}, v_{2} \bar{1})) - L((u_{1} \bar{1}, v_1\bar{1}))| \\
& \leq \ep D \min\{\la_{u_3},\la_{u_2},\ga_{v_3},\ga_{v_2}\} + \ep D \min\{\la_{u_2},\la_{u_1},\ga_{v_2},\ga_{v_1}\} \\
& \leq 2\ep e^\ep D \min\{\la_{u_3},\la_{u_1},\ga_{v_3},\ga_{v_1}\} \\
& < 4\ep D \min\{\la_{u_3},\la_{u_1},\ga_{v_3},\ga_{v_1}\}.
\end{split}
\]
Also, $\la_{u_1}/\la_{u_3},\ga_{v_1}/\ga_{v_3} \in (e^{-2\ep},e^{2\ep})$. Again, $O^\la_{u_1}= O^\la_{u_3}$ and $O^\ga_{v_1}= O^\ga_{v_3}$ are trivial; therefore $u_1 \times v_1$ and $u_3 \times v_3$ are $4\ep$-relatively close.
\end{proof}


\begin{lemma}\label{prop1}
For any $\de>0$, $r\in(0,1)$, there exists $K\in\mathbb{N}$ with the following property: For any words $s,t$ with $\la_s / \ga_t \in [r, r^{-1}]$, there are words $\al,\beta$ such that $|\al|,|\beta|\leq K$, $O^\la_{s\al}=O^\ga_{t\beta}$ and $\la_{s\al} / \ga_{t\beta} \in (e^{-\de},e^\de)$.
\end{lemma}

For the proof, we are going to use a result about convolutions of distributions. A set of real numbers is called \emph{$\sigma$-arithmetic} if all numbers in the set are integer multiples of $\sigma$ and $\sigma$ is the biggest number with this property. A distribution $F$ on $\mathbb{R}$ is called \emph{$\sigma$-arithmetic} if its support is $\sigma$-arithmetic. For $n\in\mathbb{N}$, $F^{n*}$ denotes the $n$-fold convolution of $F$ with itself.

\begin{lemma}[Feller]\label{Feller1}
Let $F$ be a distribution in $\mathbb{R}$ and $\Sigma$ be the 
union of the supports of $F,F^{2*},F^{3*},\ldots$ Assume $F$ is not concentrated on a half axis. Then,  $\Sigma$ is dense in $\mathbb{R}$ if $F$ is non-arithmetic, and $\Sigma=\{0,\pm\sigma,\pm 2\sigma,\ldots\}$ if $F$ is $\sigma$-arithmetic.
\end{lemma}
\begin{proof}
See \cite{Feller}, Vol II, Lemma V.4.2.
\end{proof}

\begin{proof}[Proof of Lemma \ref{prop1}]
First we observe that without loss of generality we may assume $O^\la_s=O^\ga_t$. Because if $O^\la_s=-O^\ga_t$, we can replace $s$ with $si_0$ where $i_0$ is any digit such that $O^\la_{i_0}=-1$ (or if no such $i_0$ exists we use $j_0$ with $O^\ga_{j_0}=-1$ and $t j_0$). Since we are claiming that the result is true for any given $r\in(0,1)$, the hypothesis of the statement will not be violated by this replacement.

We are first going to consider the case when the set 
\[
S=\{\log\la_1,\ldots,\log\la_A,-\log\ga_1,\ldots,-\log\ga_B\}
\]
is non-arithmetic. Let $R=-\log r$. The condition on $s$ and $t$ imply that $\log\la_s - \log\ga_t \in[-R,R]$. Consider a distribution $F$ whose support is the set 
\[
\begin{split}
2S&=\{2\log\la_1,\ldots,2\log\la_A,-2\log\ga_1,\ldots,-2\log\ga_B\}\\
  &=\{\log\la_{11},\ldots,\log\la_{AA},-\log\ga_{11},\ldots,-\log\ga_{BB}\}.
\end{split}
\]
Clearly, $2S$ arithmetic iff $S$ is arithmetic. If $F$ is a distribution with point masses at the elements of $2S$, by Lemma \ref{Feller1}, given $\delta>0$, there is a $K^\prime$ such that the unions of the supports of $F,F^{2*},\ldots,F^{K^\prime *}$ form a $\delta$-net in $[-R,R]$. This means that given any number $z\in[-R,R]$ there exist nonnegative integers $m_{11},\ldots,m_{AA},n_{11},\ldots,n_{BB}$ such that
\[
z + m_{11} \log\la_{11} + \cdots + m_{AA} \log\la_{AA} + n_{11}\log\ga_{11} +\cdots + n_{BB}\log\ga_{BB} \in (-\delta,\delta)
\]
and $|m_{11} +\cdots +m_{AA} +n_{11}+\cdots+n_{BB}|\leq K^\prime$. This implies that there are words $\al,\beta$ with $|\al|,|\beta|\leq 2K^\prime=:K$ satisfying 
\[
z+\log \la_\al - \log \ga_\beta \in (-\delta,\delta).
\]
We also have $O^\la_\al = O^\ga_\beta=1$ since each digit appears in pairs in these words. Taking $z=\log\la_s - \log\ga_t$ proves the result.

Now we turn to the case when $S$ is $\sigma$-arithmetic for some $\sigma>0$. Then $\log\la_s - \log\ga_t = k\sigma$ for some integer $k$ with $|k\sigma|\leq R$. Choose $\al_0,\beta_0$ such that $\la_{\al_0}=\ga_{\beta_0}$ and $O^\la_{\al_0}=-O^\ga_{\beta_0}$ if such $\al_0,\beta_0$ exist (this choice is independent of $s,t,r$ and $\delta$). Now, using the $\sigma$-arithmeticity of $S$ and arguments similar to those in the previous case (but using $S$ instead of $2S$), we can find $\al$ and $\beta$ of bounded length such that $\la_{s\al}=\ga_{t\beta}$. If $O^\la_{s\al}=O^\ga_{t\beta}$ then we are done. If not, then the words $s\al\al_0$ and $t\beta\beta_0$ will satisfy the requirements. Note that if no $\alpha_0,\beta_0$ as above exist, then we necessarily have $O^\la_{s\al}=O^\ga_{t\beta}$.
\end{proof}

\begin{corollary}\label{prop1_corollary}
For any $\de>0$, $r\in(0,1)$ there exists $K\in\mathbb{N}$ with the following property: For any words $s_1,s_2$ with $\la_{s_1} / \la_{s_2} \in [r, r^{-1}]$, there are words $\al_1,\al_2$ such that $|\al_1|,|\al_2|\leq K$, $O^\la_{s_1\al_1}=O^\la_{s_2\al_2}$ and $\la_{s_1\al_1} / \la_{s_2\al_2} \in (e^{-\de},e^\de)$.
\end{corollary}
\begin{proof}
The proof uses the same idea as above, with 
\[
S=\{\log\la_1,\ldots,\log\la_A,-\log\la_1,\ldots,-\log\la_A\}.
\]
\end{proof}

\begin{lemma}\label{cl4}
Suppose that for any $\epsilon > 0$ and $\delta > 0$ there exist two distinct $\epsilon$-relatively close $\delta$-squares. Then, given any $\epsilon >0$, $\delta > 0$ and $N \in \mathbb{N}$, there exist $N$ distinct pairwise $\epsilon$-relatively close $\delta$-squares.
\end{lemma}

\begin{proof}
The result is true for $N$=2 by assumption. Assume it is true for some $N$; we will prove that it is then true for $2N$.
Given $\epsilon$ and $\delta$, choose $0< \ep_1 < \ep/4$ and $0<\de_1 <\de$. Find $N$ distinct $\ep_1$-relatively close $\de_1$-squares $u_1 \times v_1, \ldots, u_N \times v_N$.

Now choose $\ep_2,\de_2>0$ small enough such that $\de_1+\de_2< \de$ with
\begin{equation}\label{choice1}
\ep_1 e^{\de_2} + 2 |e^{\de_2} -1| e^{\de_2 + \de_1 + \ep_1} \la_{u_1}^{-1} < \ep/4
\end{equation}

and
\begin{equation}\label{choice2}
 e^{\de_1} \la_{u_1}^{-1} e^{\ep_1} (\ep_2 + 2 e^{\ep_2 +\delta_2} |e^\ep_2-1|)< \ep/4.
\end{equation}

From our assumption that the statement is true for $N=2$, it follows that there exist $\ep_2$-relatively close $\de_2$-squares $s_1\times t_1$ and $s_2 \times t_2$.

Then, $s_1 u_i \times t_1 v_i$, $i=1,\ldots,N$ are mutually $\epsilon/4$-relatively close $\delta$-squares by Lemma~\ref{cl1} and (\ref{choice1}). Same holds with $s_2$ and $t_2$. And for any $i$, $s_1 u_i \times t_1 v_i$ and $s_2 u_i \times t_2 v_i$ are $\epsilon/4$-relatively close $\delta$-squares by Lemma~\ref{cl2} and (\ref{choice2}). Then, the $2N$ $\delta$-squares $s_i u_j \times t_i v_j$ are distinct and pairwise $\epsilon$-relatively close by Lemma~\ref{cl3} (of course we can assume $\ep<\log 2$).
\end{proof}

\bigskip

Recall that we use $d_\la,d_\ga$ to denote the similarity dimensions of the systems $\{F_i\}_1^A$ and $\{G_j\}_1^B$. Now we define a probability measure $\mu_\la$ on $\Omega_\la$ by setting $\mu_\la([u]) = \la_u^{d_\la}$ for a word $u$. Define $\mu_\ga$ similarly on $\Omega_\ga$. Note that these are well-defined by the definitions of $d_\la$ and $d_\ga$. Then $\mu = \mu_\la \times \mu_\ga$ is a probability measure on $\Omega$ satisfying  

\begin{equation}\label{measuredef}
\mu([u \times v]) = \la_u^{d_\la} \ga_v^{d_\ga}.
\end{equation}

\begin{lemma}\label{cl6}
Given any $\de>0$ and words $u$ and $v$, for $\mu$-almost all $(\om,\tau)\in\Om$, there is sequence of prefixes $\{s_j\}$ to $\om$ and $\{t_j\}$ to $\tau$ such that $\la_{s_j} / \ga_{t_j} \in (e^{-\de},e^\de)$, $O^\la_{s_j}=O^\ga_{t_j}$ and $s_j u$, $t_j v$ are prefixes of $\om$, $\tau$, respectively, for each $j$.
\end{lemma}
\begin{proof}
Recall that $r_{min}=\min\{\la_1,\ldots,\la_A,\ga_1,\ldots,\ga_B\}$. Let $R=-\log r_{min}$. Find the $K$ corresponding to $\de$ and $r_{min}$ as in Lemma \ref{prop1}. There are finitely many cylinders of the form $[s \times t]$ where $|s|,|t|\leq K$. Let $c^\prime$ be the minimum of the $\mu$-measure of such cylinders, and let $c= c^\prime \la_u^{d_\la} \ga_v^{d_\ga}$.

Observe that, any cylinder set $[s \times t]$ can be written as a disjoint union of cylinders $[s^\prime \times t^\prime]$ such that $\la_{s^\prime} / \ga_{t^\prime} \in [r_{min},1) \subset [r_{min},1/r_{min}]$ (in particular, $s^\prime \times t^\prime$ are $R$-squares). To see this, without loss of generality assume $\la_s \leq \ga_t$. Then take all shortest possible words $t^\prime$ such that $t$ is a prefix of $t^\prime$ and $\ga_{t^\prime} < \la_s$. Then the union of all subcylinders of the form $[s \times t^\prime]$ satisfy the requirements.

Now let $A\subset \Om$ be the set of points that do not have an address satisfying the statement of the lemma. For a natural number $n$, let $A_n$ be the set of points $(\om,\tau)$ of $\Om$ where $\om$ and $\tau$ have prefixes $su$, $tv$ with the desired property only for some $s,t$ satisfying $|s|,|t| \leq n$. Clearly, $A=\cup_n A_n$, therefore it suffices to show $\mu(A_n)=0$ for all $n$. Fix an $n$ now.

Write $\Om$ as a union of disjoint $R$-squares and call this collection of cylinders $S_1$. Consider a cylinder $[s \times t]$ from $S_1$. By Lemma \ref{Feller1} it has a subcylinder $[s\al \times t\beta]$ such that $|\al|,|\beta| \leq K$ and $\la_{s\al} / \ga_{t\beta} \in (e^{-\de},e^\de)$. Note that
\[
\mu([s\al u \times t\beta v]) \geq c \mu([s \times t]).
\]
Now write $[s \times t] \backslash [s\al u \times t\beta v]$ as a disjoint union of $R$-squares. Let $S_2$ be the collection of $R$-squares obtained by repeating this procedure on each $s \times t$ from $S_1$. Continue this way to construct $S_m$ for each $m$. If $\Om_m$ denotes the union of the cylinders in $S_m$, then $\{\Om_m\}$ is a decreasing sequence and
\[
\mu(\Om_{m+1}) \leq (1-c)\mu(\Om_m).
\]
Therefore $\mu(\cap \Om_m)=0$. But $A_n \subset \cap \Om_m$, therefore $\mu(A_n)=0$.
\end{proof}

\bigskip
Now let $\nu$ be the projection of $\mu$ onto $C_\la + C_\ga$ via $L$ (see (\ref{def_L}) for the definition of $L$).

\begin{proposition}\label{cl7}
If the conditions of Lemma~\ref{cl4} hold, then at $\nu$-almost all points $a$, the $(d_\la +d_\ga)$ dimensional upper density of $\nu$ is infinite, i.e. 
\[
\limsup_{r \to 0} \frac{\nu(B(a,r))}{r^{d_\ga + d_\la}} = \infty.
\]
\end{proposition}

\begin{proof}
We will prove that, given any $N \in \mathbb{N}$, the upper density is bigger than $bN$ at $\nu$-almost all points, where $b$ is a constant independent of $N$.

Given $N$, find $N$ distinct, pairwise $1/2$-relatively close $1/2$-squares $u_1 \times v_1$, \ldots, $u_N \times v_N$ (possible by Lemma~\ref{cl4}). If necessary, instead of $1/2$ we can use a small enough number so that the corresponding cylinder sets are disjoint. But for the purposes of our computation they will be regarded as $1/2$-relatively closed $1/2$-squares, so without loss of generality we may assume they are disjoint. Choose $0<\de < 1/2$ such that
\[
\frac{1}{2} e^\de + 2  |e^{\de} -1|e^{\de+1} \la_{u_1}^{-1} < 1.
\]
Let $a=L(\om,\tau)$ be a point such that for $(\om,\tau)$, $\de$ and $u_1,v_1$, there exist $s_j$ and $t_j$ as in Lemma~\ref{cl6}. By Lemma~\ref{cl6}, $\mu$-almost all $(\om,\tau)$ have this property. Then, the choice of $\delta$ as above implies that $s_j u_1 \times t_j v_1, \ldots,s_j u_N \times t_j v_N$ are pairwise 1-relatively close 1-squares by Lemma~\ref{cl1}, and these words represent disjoint cylinder sets. The image in $\mathbb{R}$ of each $[s_j u_i \times t_j v_i]$ under $L$ has diameter no more than $2De \la_{s_j u_i} \leq 2De^2 \la_{s_j u_1}$ and by 1-relative closeness, the distance from each $L(s_j u_i \bar{1},t_j v_i \bar{1})$ to $L(s_j u_1 \bar{1}, t_j v_1 \bar{1})$ is no more than
\[
D\min\{ \la_{s_j u_1}, \la_{s_j u_i}, \ga_{t_j v_1}, \ga_{t_j v_i} \} \leq e D\la_{s_j u_1}.
\]
Set $r_j=\la_{s_j u_1}$. The argument above shows that there is an absolute constant $h$ such that $B(a,hr_j)$ contains each $L([s_j u_i \times t_j v_i])$, $i=1,\ldots,N$. Observe that

\[
\begin{split}
\nu(L([s_j u_i \times t_j v_i])) &\geq \la_{s_j u_i}^{d_\la} \ga_{t_j v_i}^{d_\ga} \geq \la_{s_j u_1}^{d_\la} \ga_{t_j v_1}^{d_\ga}  \left(\frac{1}{\sqrt{e}}\right)^{d_\la + d_\ga} \\
& \geq \la_{s_j u_1}^{d_\la + d_\ga}  \left(\frac{1}{e}\right)^{d_\ga}\left(\frac{1}{\sqrt{e}}\right)^{d_\la + d_\ga} \\
& \geq r_j ^{d_\la + d_\ga} b^\prime \qquad i=1,\ldots,N 
\end{split}
\]
for some constant $b^\prime$ depending on $C_\la$ and $C_\ga$ only. Then, since the $s_j u_i \times t_j v_i$, $i=1,\ldots,N$ are disjoint in $\Om$,
\[
\frac{\nu(B(a,hr_j))}{(hr_j)^{d_\la + d_\ga}} \geq N \frac{r_j ^{d_\la + d_\ga} b^\prime}{(hr_j)^{d_\la + d_\ga}} \geq N \frac{b^\prime}{h^{d_\la + d_\ga}} = N b
\]
where $b$ is independent of $N$. Letting $j\to \infty$ (hence $r_j \to 0$), the claim is proved.
\end{proof}

This will prove one direction of the main theorem. Observe that "for every $\ep,\de$ there exist two distinct $\ep$-relatively close $\de$-squares" is equivalent to "for every $\ep^\prime$ there exist two distinct $\ep^\prime$-relatively close $\ep^\prime$-squares". To see this, one can choose $\ep^\prime=\min\{\ep,\de\}$. We state this result as a corollary:

\begin{corollary}
If for each $\ep>0$ there exist two distinct $\ep$-relatively close $\ep$-squares, then $\Ha^{d_\la+d_\ga}(C_\la+ C_\ga)=0$.
\end{corollary}
\begin{proof}
This is an easy consequence of the result above and Rogers-Taylor density theorem (for details, see \cite{Fal2}, Prop. 2.2) once we establish the fact that the restriction of $\Ha^{d_\la + d_\ga}$ to $C_\la + C_\ga$ is absolutely continuous with respect to $\nu$:
Suppose $H \subset C_\la + C_\ga$ with $\nu(H)=0$. Then (by the definition of $\nu$) $\mu(G)=0$ where $G=L^{-1}(H)$. Given any $\ep>0,\de>0$, $G$ can be covered with cylinders $u \times v$ such that $\sum \mu([u \times v]) = \sum \la_u^{d_\la} \ga_v^{d_\ga} < \ep$. By splitting each $[u \times v]$ into further subcylinders if necessary, we may assume that each $\la_u,\ga_v$ is smaller than $\de/(2D)$. Moreover, as described in Lemma~\ref{cl6}, we may also assume that $u \times v$ are $R$-squares for some absolute constant $R$, so that $\la_u \approx \ga_v$ with absolute comparison constants; implying $\la_u^{d_\ga} \ga_v^{d_\ga} \approx \la_u^{d_\la + d_\ga}$ with absolute constants. If $I_{uv}$ is the closed convex hull of $L([u \times v])$ then $|I_{uv}|< \de$ and $|I_{uv}| \approx \la_u$, again with absolute comaprison constants. Then $H\subset \bigcup I_{uv}$ and by the arguments above, $\sum |I_{uv}|^{d_\la+d_\ga} \approx \sum \la_u^{d_\la + d_\ga} \approx \sum \la_u^{d_\la} \ga_v^{d_\ga}< \ep$. Therefore we conclude that $\Ha^{d_\la + d_\ga}(H)=0$.
\end{proof}

\bigskip

Now we will prove the converse: 

\bigskip

\begin{lemma}\label{converse}
Assume that for some $\ep>0$, there does not exist two distinct $\ep$-relatively close $\ep$-squares. Then $\Ha^{d_\la +d_\ga}(C_\la + C_\ga)>0$.
\end{lemma}
\begin{proof}
Let $a\in C_\la + C_\ga$ and $r>0$. Define $r_{min}$ as in (\ref{def_rmin}). We can decompose $\Om$ into disjoint cylinders $[u \times v]$ such that 
\[
 r r_{min} \leq \la_u,\ga_v < r.
\]
In particular, such $u \times v$ are $|\log r_{min}|$-squares. Let $Y_r$ be the cylinders in this collection whose projections under $L$ intersect $B(a,r)$.

Using Lemma \ref{prop1}, find a $K\in\mathbb{N}$ corresponding to $r_{min}$ and $\ep$. Then each $[u \times v] \in Y_r$ has an $\ep$-square subcylinder $[u\al \times v\beta]$ such that $|\al|,|\beta| \leq K$. Then clearly there is an $\eta>0$ (independent of $a,u,v$ and $r$) such that

\[
\eta r \leq \la_{u\al},\ga_{v\beta} < r.
\]
(indeed we can take $\eta = {r_{min}}^{K+1}$). Let
\[
M = \left[\frac{\log(\eta r) -\log r}{\log \ep}\right] + 1 = \left[\frac{\log\eta}{\log \ep}\right] + 1.
\]

Then we have
\[
  e^{-M\ep}r \leq \eta r < e^{-(M-1)\ep} r < \cdots < e^{-2\ep} r < e^{-\ep} r < r
\]
and observe that if $t_1,t_2$ are two numbers that lie between the same two consecutive terms above, then $t_1/t_2 \in (e^{-\ep},e^\ep)$. So if we have $M+1$ numbers $t_1,\ldots,t_{M+1}$ all of which lie in the interval $(\eta r, r)$, then at least two of them satisfy $t_i/t_j \in (e^{-\ep},e^\ep)$.

How many elements can $Y_r$ have? Note that if $[u \times v] \in Y_r$ then $L(u\bar{1},v\bar{1}) \in (a-2r,a+2r)$. Then we must have

\[
\#Y_r \leq 5(M+1)^2 \frac{4r}{\eta r \ep} = (M+1)^2 \frac{20}{\eta \ep}
\]
since otherwise more than $5(M+1)^2$ of the $\ep$-subsquares $u\al \times v\beta$ will mutually satisfy condition (iii) of $\ep$-relative closeness (see Definition \ref{def2}). But then, among these, at least five will mutually satisfy the first condition as well, and two of them will certainly satisfy condition (ii) which means there exist $\ep$-relatively close $\ep$-squares, contradiction.

Note that if $[u \times v] \in Y_r$ then $\nu(L([u \times v]) = \la_u^{d_\la} \ga_v^{d_\ga} \leq r^{d_\la + d_\ga}$. Therefore

\[
\frac{\nu(B(a,r))}{r^{d_\la + d_\ga}} \leq \frac{\#Y_r r^{d_\la + d_\ga}}{r^{d_\la + d_\ga}} = \#Y_r \leq (M+1)^2 \frac{20}{\ep \eta}
\]
and the result follows again from Rogers-Taylor density theorem.
\end{proof}

\noindent {\bf Remark.} It is not difficult to verify that for every $\ep>0$ there exist two distinct $\ep$-relatively close $\ep$-squares if and only if for every $\ep>0$ there exist two distinct $\ep$-relatively close 1-squares. One direction is trivial. For the other direction, given $\ep_0>0$, find two $\ep$-relatively close 1-squares $u_i \times v_i$ where $\ep$ is much smaller than $\ep_0$. Then by Corollary~\ref{prop1_corollary} there are words $\al,\beta$ with lengths bounded by a constant $K=K(\ep_0)$ such that $u_i \al \times v_i \beta$ are $\ep_0$-squares. If $\ep$ is chosen small enough, by Lemma~\ref{cl1} $u_i \al \times v_i \beta$ will also be $\ep_0$-relatively close.

\subsection{Applications}

{\bf Sums of homogeneous sets.} One relatively simple case when it is easy to check the $\ep$-relative closeness condition of Theorem \ref{Main1} is when we have two homogeneous Cantor sets $C_\la$ and $C_\ga$ of the same diameter (recall that homogeneous means all the maps have the same contraction rate).

We will give a proof for the case when all the maps in either iterated function systems are orientation-preserving, that is, $O^\la_i=O^\ga_j=1$ for all $i,j$. If there are orientation-reversing maps, the situation is more involved; however, using the identification of $C_\la+C_\ga$ as the projection of $C_\la \times C_\ga$ on $y=x$, it can be argued geometrically that the result holds in that case, too.

We can without loss of generality assume that $C_\la$ and $C_\ga$ have diameter $1$ and their convex hulls are the interval $[0,1]$:

\begin{theorem}
Let $C_\la$ and $C_\ga$ be two homogeneous Cantor sets in the real line with closed convex hull $[0,1]$. Assume that the maps of the iterated function systems defining these sets are orientation preserving. If $d_\la$ and $d_\ga$ are the similarity dimensions of $C_\la$ and $C_\ga$, respectively, then $\Ha^{d_\la+d_\ga}(C_\la+C_\ga)=0$.
\end{theorem}

\begin{figure}[!h]

{\begin{center} \scalebox{.6}{\includegraphics{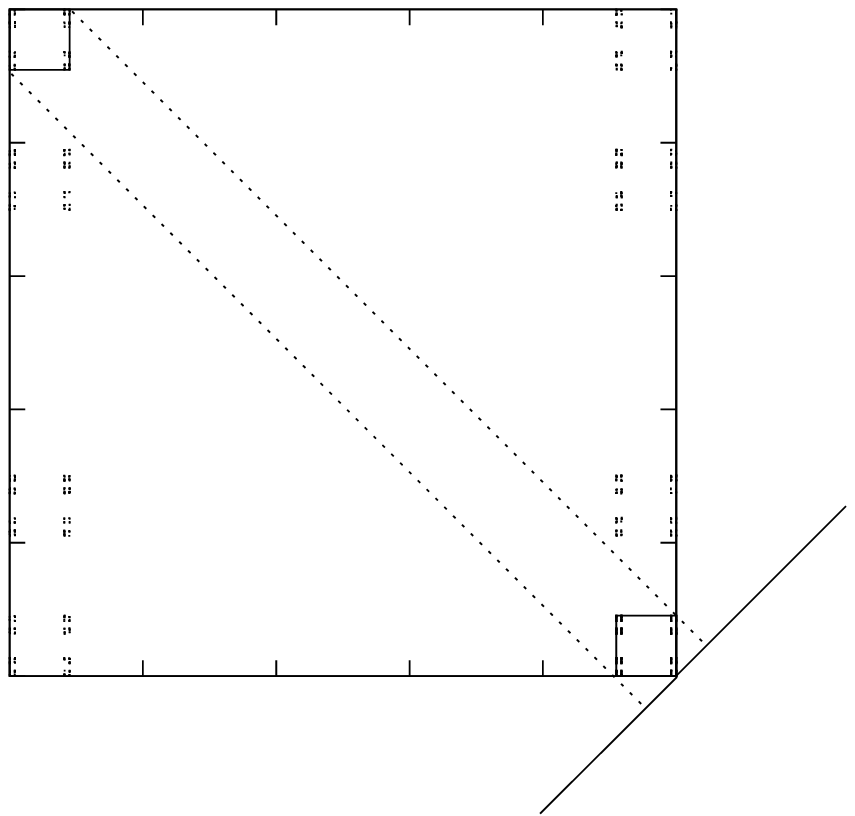}}\end{center}}
\caption{}\label{CornerOverlap}
\end{figure}

\begin{proof}
Let $\la$ and $\ga$ be the contraction rates of the maps for these two sets. Then there exist maps $F_{i_L},F_{i_R},G_{j_L},G_{j_R}$ given by:
\[
\begin{aligned}
F_{i_L}(x) &= \la x, \\
F_{i_R}(x) &= \la x + 1-\la, \\
G_{j_L}(x) &= \ga x, \\
G_{j_R}(x) &= \ga x + 1-\ga.
\end{aligned}
\]
Given $\ep>0$, we can find $n,m$ such that $\la^n/\ga^m\in(e^{-\ep},e^\ep)$ and $|1-\frac{\la^n}{\ga^m}|<\ep$. Without loss of generality, we may assume $\ga^m\leq\la^n$. Then, it is easy to check that
\[
\begin{split}
&|L(i_R^n\bar{1},j_L^m\bar{1}) -L(i_L^n\bar{1},j_R^m\bar{1})|\\
&= |\Pi_\la(i_R^n\bar{1})-\Pi_\la(i_L^n\bar{1})+ \Pi_\ga(j_L^m\bar{1})-\Pi_\ga(j_R^m\bar{1})|\\
&=|\la^n \Pi_\la(\bar{1}) + 1 -\la^n - \la^n\Pi_\la(\bar{1}) + \ga^m \Pi_\ga(\bar{1}) - (\ga^m\Pi_\ga(\bar{1}) + 1 -\ga^m)|\\
&=|\la^n - \ga^m| = \ga^m \left| 1- \frac{\la^n}{\ga^m}\right| \leq \ep \min\{\ga^m,\la^n\}
\end{split}
\]
and this means that $i_R^n \times j_L^m$ and $i_L^n \times j_R^m$ are $\ep$-relatively close $\ep$-squares. Geometrically, they correspond to the affine copies of $C_\la \times C_\ga$ at the upper left and lower right corners of $[0,1]\times[0,1]$ (see Fig. \ref{CornerOverlap}).
\end{proof}

\bigskip

{\bf Orthogonal projections.} Another application is related to the orthogonal projections of the set $C_\la \times C_\ga$ onto lines. The projection onto the line with direction vector $(1,\eta)$ is, up to scaling, the set $C_\la+\eta C_\ga$, and this is an arithmetic sum of two Cantor sets.

The result we are going to prove will be similar to an earlier result by Y.~Peres, K.~Simon and B.~Solomyak \cite{SolScan}. Let $\Pi_\theta$ be the orthogonal projection from the plane to the line $l_\theta$ through the origin making angle $\theta$ with the positive $x$-axis. For any set $E$ in the plane we can define the set $\ip$ of intersection parameters as
\[
\ip=\ip(E) = \{\theta \mid \Pi_\theta \text{ is not one-to one on $E$}\}.
\]
In other words, $\ip$ is the set of all directions that can be obtained by joining pairs of points in $E$. The following result has been established:

\begin{numarasizteo}[Peres, Simon, Solomyak]\label{SolomIP}
Let $E\subset \mathbb{R}^2$ be a self-similar set without rotation of similarity dimension $s \in (0,1)$ that is not contained in a line. Then, $\Ha^s(\Pi_\theta E)=0$ for Lebesgue-almost all $\theta\in\ip$. 
\end{numarasizteo}
Here, ``without rotation'' means that the similarity maps in the iterated function system defining $E$ are non-rotating.
Our result is as follows:
\begin{theorem}\label{sumproj}
Let $C_\la$, $C_\ga$ be affine Cantor sets in the real line with similarity dimensions $d_\la,d_\ga$, respectively. Let $\ip$ be the intersection parameter set for $C_\la \times C_\ga \subset \mathbb{R}^2$. Then, $\Ha^{d_\la+d_\ga}(\Pi_\theta (C_\la\times C_\ga))=0$ for Lebesgue almost all $\theta\in\ip$.
\end{theorem}

Of course, the result above has content if $\mathcal{L}(\ip)>0$. Note that in the special case $\Ha^{d_\la}(C_\la)>0$ and $\Ha^{d_\ga}(C_\ga)>0$ with $d_\la +d_\ga=1$ and $d_\la,d_\ga>0$, then $C_\la \times C_\ga$ is an irregular 1-set in the plane hence projects to Lebesgue null sets in Lebesgue-almost all directions, by the well-known result of Besicovitch \cite{Fal1}. Now let us investigate what happens if the $\ip$ set is small. We observe that the situation $d_\la + d_\ga = 1/2$ is a borderline for the $\ip$ set in the following sense: If $d_\la + d_\ga <1/2$, then the self-difference set $(C_\la \times C_\ga) - (C_\la \times C_\ga)$ whose circular projection can be identified with the $\ip$ set, has Hausdorff dimension less than one, which implies that $\mathcal{L}(\ip)=0$. Our second result concerns this case:

\begin{theorem}\label{smallsum}
Let $C_\la$, $C_\ga$ be affine Cantor sets in the real line satisfying strong separation condition, with similarity dimensions $d_\la,d_\ga$, respectively. Assume $d_\la + d_\ga <1/2$. Then, $\Ha^{d_\la+d_\ga}(\Pi_\theta (C_\la\times C_\ga))>0$ for Lebesgue almost all $\theta\in[0,\pi]$.
\end{theorem}

\bigskip

We first prove Theorem \ref{sumproj}: By symmetry, it suffices to prove the result for $\theta\in(0,\frac{\pi}{2})$. Recall that, up to scaling, $\Pi_\theta (C_\la\times C_\ga)$ is equal to $C_\la + \eta C_\ga$ where $\eta=\tan\theta$. We will say $\eta\in\tan\ip$ if $\eta=\tan\theta$ for some $\theta\in\ip \cap (0,\frac{\pi}{2})$. The theorem will be proved if we can show that $\Ha^{d_\la+d_\ga}(C_\la + \eta C_\ga)=0$ for Lebesgue almost all $\eta\in\tan\ip$.

Given $\ep>0$, let $P_\ep$ be the set of $\eta$ such that there does not exist two distinct $\ep$-relatively close $\ep$-squares for $C_\la$ and $\eta C_\ga$. We are going to prove that \mbox{$\mathcal{L}(\tan\ip\cap P_\ep)=0$} for all $\ep>0$. This implies that for Lebesgue almost all $\eta\in\tan\ip$ there exist $\ep$-relatively close $\ep$-squares and this proves the theorem. To prove that \mbox{$\mathcal{L}(\tan\ip\cap P_\ep)=0$}, we will show that the Lebesgue density of $P_\ep$ is less than one at each $\eta\in\tan\ip$.

Now we introduce some additional notation. For sequences (or words) $\om_1,\om_2$ of the same symbol space, we will denote by $\om_1 \wedge \om_2$ the longest common prefix of $\om_1$ and $\om_2$. For $(\om,\tau)\in\Om$, define
\begin{equation}\label{Leta}
L_\eta (\om,\tau) = \Pi_\la(\om) + \eta \Pi_\ga(\tau)
\end{equation}
(so the map $L$ we have used before is $L_1$ with this notation). A simple observation is that
\begin{equation}\label{wedges}
|L_\eta(\om_1,\tau_1) - L_\eta (\om_2,\tau_2)|\leq D_\la \la_{\om_1 \wedge \om_2} + \eta D_\ga  \ga_{\tau_1 \wedge \tau_2}.
\end{equation}

Now fix an $\ep>0$. Let $\eta_0 \in\tan\ip$. Note that $\eta_0\neq0$ as we are considering directions in $(0,\frac{\pi}{2})$. Then, there exist $(\om_1,\tau_1),(\om_2,\tau_2)\in\Om$ such that $\Pi(\om_1,\tau_1)$ and $\Pi(\om_2,\tau_2)$ are distinct points in the plane and

\begin{equation}\label{Lequality}
L_{\eta_0}(\om_1,\tau_1) - L_{\eta_0}(\om_2,\tau_2) = 0.
\end{equation}

Since we are ruling out vertical projections (by assuming $\eta_0 \ne 0$), we can write
\begin{equation}\label{c1}
|\Pi_\ga(\tau_1)-\Pi_\ga(\tau_2)| =: c_1 \neq 0.
\end{equation}

Let $r_{min}$ be as in (\ref{def_rmin}). Choose a number $\rho > 0$ small enough satisfying 
\begin{equation}\label{rho}
c_1 /2 > 2 D \frac{\rho}{\eta_0} >0
\end{equation}
where $D=\max\{D_\la,D_\ga\}$ as before. We can find prefixes $\tilde{s}_i,\tilde{t}_i$ of $\om_i,\tau_i$, respectively, $i=1,2$, such that
\begin{equation}\label{comparable}
\begin{split}
\rho r_{min} &\leq \la_{\tilde{s}_i} < \rho,\ i=1,2\\
\rho r_{min} &\leq \eta_0 \ga_{\tilde{t}_i} < \rho,\ i=1,2.
\end{split}
\end{equation}
Applying Corollary \ref{prop1_corollary} to each of the pairs $(\tilde{s}_1,\tilde{s}_2)$ and $(\tilde{t}_1,\tilde{t}_2)$ with $r_{min}$ and $\ep/4$, we can find words $\al_i,\beta_i,\ i=1,2$ of length at most $K=K(\ep,r_{min})$ such that 
\begin{equation}\label{ep3-1}
\begin{split}
\frac{\la_{\tilde{s}_1 \al_1}}{\la_{\tilde{s}_2 \al_2}} \in (e^{-\ep/4},e^{\ep/4})\ & \text{    and    }\  O^\la_{\tilde{s}_1 \al_1}=O^\la_{\tilde{s}_2 \al_2}, \\
\frac{\ga_{\tilde{t}_1 \beta_1}}{\ga_{\tilde{t}_2 \beta_2}} \in (e^{-\ep/4},e^{\ep/4})\ & \text{    and    }\  O^\ga_{\tilde{t}_1 \beta_1}=O^\ga_{\tilde{t}_2 \beta_2}.
\end{split}
\end{equation}
Note that in this case we have 
\begin{equation}
\frac{\la_{\tilde{s}_i \al_i}}{\eta_0 \ga_{\tilde{t}_i \beta_i}} \in (e^{-\Delta},e^\Delta)
\end{equation}
where $\Delta = \Delta(\ep,r_{min},\eta_0)= K \log r_{min} + \ep/4 + |\log \eta_0|$. Applying Lemma \ref{prop1} to $\tilde{s}_1\al_1, \tilde{t}_1\beta_1$ with $r_{min}$, $\Delta$ and $\ep/3$, we can find words $\al,\beta$ of length at most $K=K(\ep,r_{min},\eta_0)$ such that 
\begin{equation}\label{ep3-2}
\frac{\la_{\tilde{s}_1 \al_1 \al}}{\eta_0 \ga_{\tilde{t}_1 \beta_1 \beta}} \in (e^{-\ep/4},e^{\ep/4}).
\end{equation}
Finally, there is an open interval $I=I(\ep,\eta_0)$ about $\eta_0$ such that
\begin{equation}\label{etabound}
\frac{\eta}{\eta_0}\in (e^{-\ep/4},e^{\ep/4}).
\end{equation}
Then, combining (\ref{ep3-1}),(\ref{ep3-2}),(\ref{etabound}), we conclude that
\begin{equation}\label{allinone}
\begin{split}
\frac{\la_{\tilde{s}_1 \al_1 \al}}{\la_{\tilde{s}_2 \al_2 \al}} \in (e^{-\ep},e^{\ep})\ & \text{    and    }\  O^\la_{\tilde{s}_1 \al_1 \al}=O^\la_{\tilde{s}_2 \al_2 \al}, \\
\frac{\ga_{\tilde{t}_1 \beta_1 \beta }}{\ga_{\tilde{t}_2 \beta_2} \beta} \in (e^{-\ep},e^{\ep})\ & \text{    and    }\  O^\ga_{\tilde{t}_1 \beta_1 \beta}=O^\ga_{\tilde{t}_2 \beta_2 \beta}, \\
\frac{D_\la \la_{\tilde{s}_i \al_i \al}}{\eta D_\ga \ga_{\tilde{t}_i \beta_i \beta}} \in (e^{-\ep},e^{\ep})\ & \text{    for    }\  i=1,2,\ \eta\in I.
\end{split}
\end{equation}
Also observe that 
\begin{equation}\label{c0}
\la_{\tilde{s}_i \al_i \al},\eta \ga_{\tilde{t}_i \beta_i \beta} \geq \rho c_0,\ i=1,2,\ \eta\in I
\end{equation}
where $c_0=c_0(\ep,r_{min},\eta_0)$ by (\ref{comparable}) and the boundedness of the lengths of $\al_i,\al,\beta_i,\beta$.

Set $s_i=\tilde{s}_i\al_i\al$ and $t_i=\tilde{t}_i\beta_i\beta$, $i=1,2$. Then (\ref{comparable}), (\ref{allinone}) and (\ref{c0}) imply that, $s_1 \times t_1$ and $s_2\times t_2$ are $\ep$-relatively close $\ep$-squares for all $\eta \in I$ satisfying
\begin{equation}
|L_\eta (s_1\bar{1},t_1\bar{1}) - L_\eta (s_2\bar{1},t_2\bar{1})| \leq \ep D c_0 \rho.
\end{equation}

Let $\Phi(\eta)=L_\eta (s_1\bar{1},t_1\bar{1}) - L_\eta (s_2\bar{1},t_2\bar{1})$. We will find an interval near $\eta_0$ where $\Phi(\eta)$ is small. We first observe that using (\ref{Leta}), (\ref{c1}), (\ref{rho}) and (\ref{comparable}), on the interval $I$ we have

\begin{equation}\label{c2}
\begin{split}
|\Phi^\prime(\eta)| &= \left| \frac{d}{d\eta}\left( \Pi_\la(s_1\bar{1})-\Pi_\la(s_2\bar{1}) + \eta \left(\Pi_\ga(t_1\bar{1}) - \Pi_\ga(t_2\bar{1})\right) \right) \right| \\
& = | \Pi_\ga(t_1\bar{1}) - \Pi_\ga(t_2\bar{1}) | \\
&\geq  | \Pi_\ga(\tau_1) - \Pi_\ga(\tau_2) | - |(\Pi_\ga(\tau_1)- \Pi_\ga(t_1\bar{1})) -(\Pi_\ga(\tau_2)- \Pi_\ga(t_2\bar{1})) | \\
&\geq c_1 - D_\ga \ga_{\tau_1 \wedge t_1} - D_\ga \ga_{\tau_2 \wedge t_2} \geq c_1 - D_\ga \ga_{\tilde{t_1}} - D_\ga \ga_{\tilde{t_2}}\\
&\geq c_1 -2 D \frac{\rho}{\eta_0} > c_1 / 2 =:c_2 >0.
\end{split}
\end{equation}
Using the triangle inequality in the other direction we can also easily get
\begin{equation}\label{c3}
|\Phi^\prime(\eta)| \leq 2c_1 =:c_3.
\end{equation}

Next we find a bound for $\Phi(\eta_0)$. Using (\ref{wedges}), (\ref{Lequality}) and (\ref{comparable}), we get
\begin{equation}\label{c4}
\begin{split}
|\Phi(\eta_0)| = & |\Phi(\eta_0)- (L_{\eta_0}(\om_1,\tau_1) - L_{\eta_0}(\om_2,\tau_2))|\\
=&| L_{\eta_0}(s_1\bar{1},t_1\bar{1})- L_{\eta_0}(\om_1,\tau_1) + (L_{\eta_0}(s_2\bar{1},t_2\bar{1}) - L_{\eta_0}(\om_2,\tau_2))|\\
\leq& D_\la \la_{s_1 \wedge \om_1} + \eta_0 D_\ga \ga_{t_1 \wedge \tau_1} + D_\la \la_{s_2 \wedge \om_2} + \eta_0 D_\ga \ga_{t_2 \wedge \tau_2} \\
\leq& D_\la \la_{\tilde{s_1}} + \eta_0 D_\ga \ga_{\tilde{t_1}} + D_\la \la_{\tilde{s_2}} + \eta_0 D_\ga \ga_{\tilde{t_2}}\\
\leq& 4 D (1+\eta_0) \rho =: c_4 \rho.
\end{split}
\end{equation}
But then, provided we are in the interval $I$, (\ref{c2}) and (\ref{c4}) imply that $\Phi(\eta^\prime)=0$ for some $\eta^\prime\in B(\eta_0,c_4\rho/c_2)$, and (\ref{c3}) implies that 
\begin{equation}\label{fismall}
|\Phi(\eta)| < \ep D c_0 \rho \ \ \text{ if \ $\eta\in B\left(\eta^\prime,\frac{\ep D c_0 \rho}{c_3}\right)$}
\end{equation}
and observe that such $\eta$ are outside the set $P_\ep$.
Since $ B\left(\eta_0,\frac{\ep D c_0 \rho}{c_3}+\frac{c_4 \rho}{c_2}\right) \subset I$ holds for $\rho$ small, by taking arbitrarily small $\rho$ and noting that the constants $c_i$ are independent of $\rho$, we see that the Lebesgue density of $P_\ep$ at $\eta_0$ is at most
\[
1 - \frac{\frac{\ep D c_0 \rho}{c_3}}{\frac{c_4 \rho}{c_2}+\frac{\ep Dc_0 \rho}{c_3}} = 1 - \frac{\frac{\ep Dc_0 }{c_3}}{\frac{c_4}{c_2}+\frac{\ep Dc_0 }{c_3}} < 1
\]
and therefore the theorem is proved. \qed

\bigskip

Now we return to Theorem \ref{smallsum}: Recall that we are now assuming strong separation, and that $d_\la + d_\ga <1/2$. Given $M>0$ and $\ep>0$, let $\Theta(M,\ep)$ be the set of angles $\theta$ for which $|\tan \theta| \leq M$ and there exist two distinct $\ep$-squares that are $\ep$-relatively close in the projection onto the line $l_\theta$. We will identify the projection with the set $C_\la + \eta C_\ga$. We will prove that $\lim_{\ep\to 0}\mathcal{L}(\Theta(M,\ep))= 0$, and since $M$ is arbitrary, this implies the desired result by Theorem \ref{Main1}.

Let $0<\rho<1$ and set $I_k = (\rho^k,\rho^{k-1}]$. Fix an $\ep>0$ small. Define

\[
S_k = \{ s \times t \mid s\times t \text{ is an $\ep$-square and } \la_s \in I_k\}.
\]
and let $P_k$ be the set of all $\ep$-square pairs $(s_1\times t_1, s_2 \times t_2)$ with $s_1\times t_1 \in S_k$ that also satisfy the first two conditions of Definition \ref{def2}. Since in this case (assuming $\ep$ is small enough) we have $s_2\times t_2 \in S_{k-1}\cup S_k \cup S_{k+1}$, it follows that $(\# P_k) \leq C_1 (\# S_k)^2$ where $C_1$ is independent of $k$ and $\ep$. To estimate $\# S_k$, observe that, if $\la_s \in I_k$ then $\mu_\la ([s]) = \la_s^{d_\la} \asymp \rho^{k d_\la}$. We can choose $\rho$ to be sufficiently close to $1$ so that if $\la_{s_1},\la_{s_2}\in I_k$ are distinct words, then neither is a prefix of the other, hence $\mu_\la([s_1] \cap [s_2])=0$. This implies that
\[
\#\{s\mid \la_s \in I_k\} \asymp \rho^{-k d_\la}.
\]
If $s \times t$ is an $\ep$-square for $C_\la + \eta C_\ga$ then $\la_s \in I_k$ implies (again, for $\ep$ small) $\eta\la_t \in I_{k-1} \cup I_k \cup I_{k+1}$, Since $|\eta|=|\tan\theta|\leq M$, by making $\rho$ even closer to $1$ if necessary and assuming $\ep$ is small, for a fixed $s$ with $\la_s\in I_k$, we have
\[
\#\{t \mid \text{$s \times t$ is an $\ep$-square for $C_\la + \eta C_\ga$}\} \asymp \rho^{-k d_\ga}
\]
and the proportionality constants in the last two equations are independent of $k$ and $\ep$. Therefore, we conclude that there exists $C_2 >0 $ independent of $k$ and $\ep$ such that
\[
\# P_k \leq C_1 (\# S_k)^2 \leq C_2 (\rho^{-k d_\la} \rho^{-k d_\ga})^2 = C_2 \rho^{-2k(d_\la + d_\ga)}.
\]

Given a pair $(s_1\times t_1, s_2 \times t_2)$ from $P_k$, let $\theta_0$ be the direction of the line segment that joins $\Pi(s_1\bar{1},t_1\bar{1})$ to $\Pi(s_2\bar{1},t_2\bar{1})$. Note that the distance between these two points in the plane is comparable to $\la_{s_1}$, hence to $\rho^k$, by strong separation. A simple geometric argument will show that, there is a constant $C_3$ (independent of $k,\ep$) such that the projections of $s_1 \times t_1$ and $s_2\times t_2$ onto $l_\theta$ cannot be $\ep$-relatively close if $|\theta - \theta_0| > C_3 \ep \rho^k$. Therefore, the set of angles for which there is a pair $(s_1\times t_1, s_2 \times t_2)$ of $\ep$-squares that are projected to $\ep$-relatively close squares where $s_1\times t_1 \in P_k$, has Lebesgue measure no more than
\[
 C_2 \rho^{-2k(d_\la + d_\ga)} C_3 \ep \rho^k = C_2 C_3 \ep (\rho^k)^{1- 2(d_\la + d_\ga)} =: C_4 \ep \rho^{k\delta}
\]
where $\delta >0$. Therefore,
\[
\mathcal{L}(\Theta(M,\ep)) \leq  \sum_{k=1}^\infty C_4 \ep \rho^{k\delta} =: C_5 \ep
\]
where $C_5$ doesn't depend on $\ep$. This proves that $\mathcal{L}(\Theta(M,\ep))\to 0$ as $\ep\to 0$, and the result follows. \qed

\bigskip

\begin{figure}[h]
\psfrag{lam}{$\la$}
\psfrag{gam}{$\ga$}
\psfrag{eq1}{$\frac{\la}{1-2\la}\frac{\ga}{1-2\ga}=1$}
\psfrag{eq2}{$d_\la+d_\ga=1$}
\psfrag{eq3}{$\frac{2\la}{1-3\la}\frac{2\ga}{1-3\ga}=1$}
\psfrag{eq4}{$d_\la+d_\ga=\frac{1}{2}$}
\scalebox{0.8}{\includegraphics{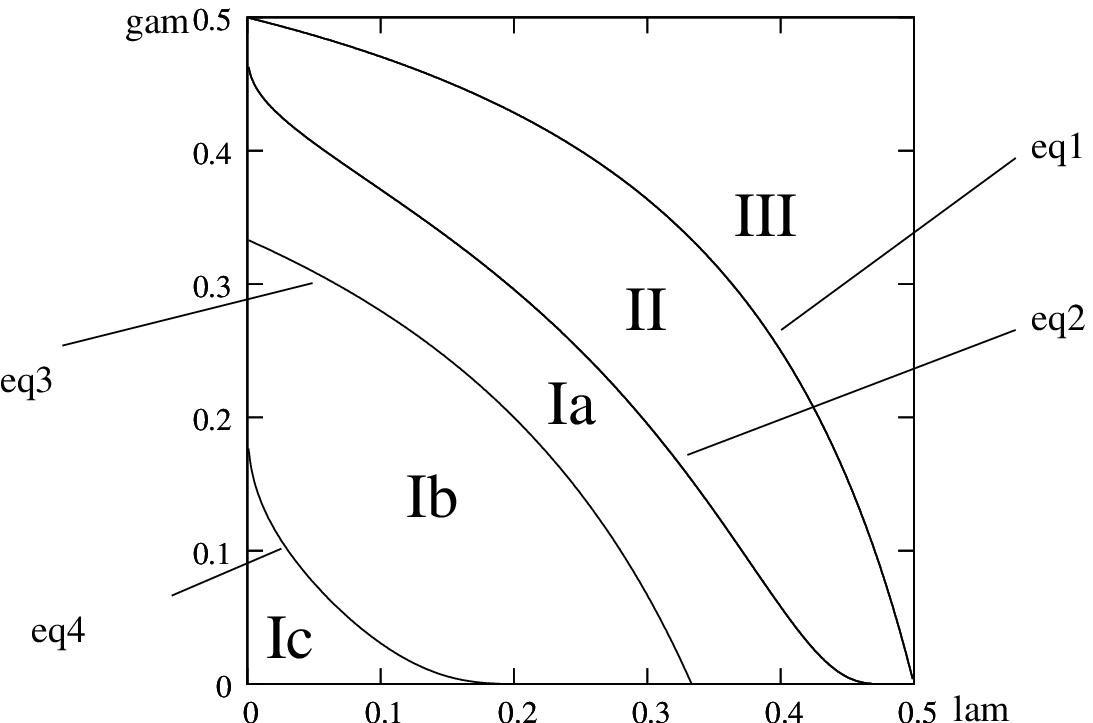}}
\caption{}\label{bolgeler}
\end{figure}

We conclude by illustrating the results above on an example that was studied by several other authors: We will take $C_\la$ and $C_\ga$ to be the middle-$(1-2\la)$ and middle-$(1-2\ga)$ cantor sets, respectively. The Newhouse thicknesses \cite{Newhouse} of these sets are $\la/(1-2\la)$ and $\ga/(1-2\ga)$. When the product of these numbers is bigger than $1$, it follows that $C_\la+C_\ga=[0,2]$ since the sets have equal diameter (see \cite{Newhouse} for the proofs of the related theorems). The set of such $(\la,\ga)$ have been indicated as region III in Figure \ref{bolgeler}. Another critical case is when $d_\la + d_\ga=1$. We have seen that if the sum of the dimensions is bigger than $1$ then typically the sum has positive Lebesgue measure. The curve $d_\la + d_\ga=1$ represents the lower boundary of region II in the figure. In this region, for almost all parameter pairs on each vertical or horizontal line we have $\mathcal{L}(C_\la + C_\ga)>0$ by Theorem~\ref{SolInt1}. Mendes and Oliveira \cite{Mendes} have shown that there are countably many ``spikes'' extending from the common boundary of regions II and III into region II, on which the sum is the full interval $[0,2]$. In region I we have $d_\la + d_\ga <1$, so clearly the sum has zero Lebesgue measure. We divide it into subregions to study the projection properties:

Let $\theta_0=-\pi/4$, the direction of the line segment joining $(0,1)$ to $(1,0)$. Consider the affine copies of $C_\la \times C_\ga$ at the upper-left and lower-right corners of $[0,1]\times[0,1]$, of size $\la^m \times \ga^n$ where $\la^m/\ga^n$ is very close to $1$ (it can be chosen arbitrarily close to $1$). The projections of these copies (up to scaling) are translates, say, $T_1(t_0)$ and $T_2(t_0)$, of $\la^m C_\la + \ga^n t_0 C_\ga$ with $t_0=\tan \theta_0=1$, and the interiors of the convex hulls of these sets overlap. As $\theta$ changes on a small interval about $\theta_0$, rescaling, we obtain a family of sets $T_1(t)$ and $T_2(t)=T_1(t)+\sigma(t)$ that are translates of $\la^m C_\la + \ga^n t C_\ga$, the relative translation $\sigma(t)$ being continuous in $t$. Thus the convex hulls of $T_1(t)$ and $T_2(t)$ overlap when $t$ is near $t_0$. If it were true that $T_1(t)-T_1(t)$ is an interval, then it would follow that $T_1(t)\cap T_2(t)$ is nonempty, which implies that the angle corresponding to $t$ is in the $\ip$ set, and this would prove that the $\ip$ set contains an interval. But

\[
T_1(t)-T_1(t)=(\la^m C_\la + t \ga^n C_\ga)-(\la^mC_\la + t \ga^n C_\ga)= \la^m (C_\la -C_\la)+t\ga^n (C_\ga-C_\ga).
\]
We want to prove that the sum on the right is a full interval. The summands $\la^m (C_\la -C_\la)$ and $t\ga^n(C_\ga-C_\ga)$ are self-similar sets with Newhouse thicknesses $2\la/(1-3\la)$ and $2\ga/(1-3\ga)$, respectively. Their diameters are $2\la^m$ and $2t\ga^n$. Since we are assuming that the ratio of these numbers is sufficiently close to $1$, by the results of Newhouse \cite{Newhouse} we can conclude that the sum is an interval provided
\[
\frac{2\la}{1-3\la}\frac{2\ga}{1-3\ga}>1.
\]
This is the inequality that defines the region Ia in Figure \ref{bolgeler}. In that region the $\ip$ set contains an interval, hence $\mathcal{L}(\ip)>0$.

Finally, regions Ib and Ic are separated by the curve $d_\la + d_\ga=1/2$. We have seen that in region Ic, Lebesgue almost all projections have positive $(d_\la + d_\ga)$-dimensional Hausdorff measure.

The situation in region Ib (together with its boundaries) is not known. It is expected that $\mathcal{L}(\ip)>0$ there. Peres, Simon and Solomyak \cite{SolScan} have investigated the case $\la=\ga=r$ (which gives a self-similar four-corner set) and found that if $1/6 < r < 1/4$ then $\ip$ contains an interval. That $\dim \ip < 1$ if $r<1/9$ is seen easily using self-similarity. The case $1/9 < r < 1/6$ is still open, as is the question on the rest of region Ib.

\bigskip
\bigskip

\noindent {\bf Acknowledgements.} The author would like to thank Boris Solomyak for his valuable comments and suggestions.

\bibliographystyle{amsplain}
\bibliography{general}

\end{document}